\documentclass[12pt]{article}
\usepackage{amsmath, amsthm, amssymb}
\numberwithin{equation}{section}
\begin{document}
\author{Ajai Choudhry\\
and\\
Arman Shamsi Zargar}
\title{A sextic   diophantine  chain \\and  a related Mordell curve }
\date{}
\maketitle

\begin{abstract}
In this paper we obtain  parametric as well as numerical solutions of the sextic   diophantine  chain $
\phi(x_1,\,y_1,\,z_1)=\phi(x_2,\,y_2,\,z_2)=\phi(x_3,\,y_3,\,z_3)=k$ where $\phi(x,\,y,\,z)$ is a sextic   form
defined by $\phi(x,\,y,\,z)$ $=x^6+y^6+z^6-2x^3y^3-2x^3z^3-2y^3z^3$ and $k$ is an integer.  Each numerical solution of such a  sextic   chain yields, in general, nine rational points on the  Mordell curve $y^2=x^3+k/4$. While all of these nine points are not independent in the group of rational points of the Mordell curve, we have constructed a parameterized family of Mordell curves of generic rank $\geq 6$ using the aforementioned parametric solution of the sextic   diophantine  chain. Similarly, the numerical solutions of the sextic   chain yield  additional examples of Mordell curves whose rank is $\geq 6$. 
\end{abstract}



\section{Introduction}
Let $\phi(x,\,y,\,z)$ be a sextic   form, with integer coefficients, in the three variables $x,\,y$ and $z$. While a limited number of diophantine  equations of the type $\phi(x_1,\,y_1,\,z_1)=\phi(x_2,\,y_2,\,z_2)$ have been solved (see \cite{Br}, \cite{Ch1}, \cite{Ch2}), until now no sextic   diophantine  chains of the type 
\begin{equation}
\phi(x_1,\,y_1,\,z_1)=\phi(x_2,\,y_2,\,z_2)=\phi(x_3,\,y_3,\,z_3), \label{sexticchn}
\end{equation}
have been published. 

In this paper we obtain  parametric as well as numerical solutions of the diophantine chain~\eqref{sexticchn} where the form $\phi(x,\,y,\,z)$  is defined by
\begin{equation}
\phi(x,\,y,\,z)=x^6+y^6+z^6-2x^3y^3-2x^3z^3-2y^3z^3. \label{sexticphi}
\end{equation}

It is interesting to observe that if we write $k=\phi(x,\,y,\,z)/4$ where $x,\,y,\,z$ are rational numbers, three rational points on the Mordell curve
\begin{equation}
y^2=x^3+k, \label{Mordelleq}
\end{equation}
are given by 
\[(xy,\;(x^3+y^3-z^3)/2),\quad (yz,\;(y^3+z^3-x^3)/2),\quad (xz,\;(x^3+z^3-y^3)/2).\]

In view of the above,  it is clear  that any solution of the diophantine  chain~\eqref{sexticchn} immediately yields 9 rational points (not necessarily distinct) on the Mordell curve~\eqref{Mordelleq} where $k=\phi(x_1,\,y_1,\,z_1)/4$. Thus, apart from their intrinsic interest,  solutions of the sextic   diophantine  chain~\eqref{sexticchn} are  also   expected to yield  Mordell curves of high rank. 

In this context it is pertinent to note that Kihara \cite{Ki} has obtained a family of Mordell curves over the field $\mathbb{Q}(t)$ of generic rank $\geq 6$ and till now, this is the best known result of this type  regarding families of Mordell curves.  The parametric solution of the sextic   diophantine  chain~\eqref{sexticchn} obtained in this paper  also yields a family of Mordell curves, defined  over the field $\mathbb{Q}(m,\,n)$,  of generic rank $\geq 6$.

\section{Sextic diophantine  chains}

In Section~2.1 we describe a general method of constructing  certain sextic   diophantine  chains. We will apply this method in Sections~2.2 and 2.3 to obtain solutions of the sextic   diophantine  chain~\eqref{sexticchn} where $\phi(x,\,y,\,z)$ is defined by \eqref{sexticphi}. We also show how more  solutions of the sextic   chain~\eqref{sexticchn} may be obtained.

\subsection{A general method of constructing sextic    diophantine  chains}
Let an arbitrary  sextic   form $\phi(x,\,y,\,z)$ in  the three variables $x,\,y,\,z$ be expressible as
\begin{equation}
\phi(x,\,y,\,z)=k_1Q^3(x,\,y,\,z)+k_2C^2(x,\,y,\,z), \label{sexticgen}
\end{equation}
where $k_1,\,k_2$ are constants while $Q(x,\,y,\,z)$ is a quadratic form and $C(x,\,y,\,z)$ is a cubic form in the three variables $x,\,y,\,z$ such that
\begin{equation}
Q(x,\,y,\,z)=Q(y,\,x,\,z) \quad {\rm and}\quad  C(x,\,y,\,z)=C(y,\,x,\,z). \label{condQC}
\end{equation}

It is clear from \eqref{sexticgen} that to construct a sextic   diophantine  chain~\eqref{sexticchn}, it suffices to construct the simultaneous diophantine  chains,
\begin{align}
Q(x_1,\,y_1,\,z_1)&=Q(x_2,\,y_2,\,z_2)=Q(x_3,\,y_3,\,z_3), \label{qdchngen}\\
C(x_1,\,y_1,\,z_1)&=C(x_2,\,y_2,\,z_2)=C(x_3,\,y_3,\,z_3).  \label{cubchngen}
\end{align}

We will first obtain a parametric solution of  the simultaneous diophantine  equations,
\begin{equation}
\begin{aligned}
Q(x_1,\,y_1,\,z_1)&=Q(x_2,\,y_2,\,z_2),\\
C(x_1,\,y_1,\,z_1)&=C(x_2,\,y_2,\,z_2), \label{qdcubeqn}
\end{aligned}
\end{equation}
together with the auxiliary equation,
\begin{equation}
x_1+y_1+hz_1=x_2+y_2+hz_2, \label{habgen}
\end{equation}
where $h$ is a rational parameter,  and then use this solution to obtain a solution of the simultaneous diophantine  chains~\eqref{qdchngen} and \eqref{cubchngen} together with  the auxiliary diophantine  chain,
\begin{equation}
x_1+y_1+hz_1=x_2+y_2+hz_2 = x_3+y_3+hz_3. \label{chnhabc}
\end{equation}

A parametric solution of the simultaneous equations~\eqref{qdcubeqn} and \eqref{habgen} may be obtained either by the general method described in \cite{Ch0} or by any other appropriate method. 

Now, let $(x_i,\,y_i,\,z_i)=(\alpha_i,\,\beta_i,\,\gamma_i),\;i=1,\,2$, be a solution of the simultaneous equations~\eqref{qdcubeqn} and  \eqref{habgen} where $\alpha_i,\,\beta_i,\,\gamma_i,\;i=1,\,2$ are given in terms of certain independent parameters. We will solve the simultaneous equations~\eqref{qdchngen}, \eqref{cubchngen} and \eqref{chnhabc} by obtaining three distinct solutions of the following three simultaneous equations,
\begin{align} 
x+y+hz&=k_1, \label{eqk1}\\
Q(x,\,y,\,z)&=k_2, \label{eqk2}\\
C(x,\,y,\,z)&=k_3, \label{eqk3} 
\end{align}
where
\begin{equation} 
\begin{aligned}
k_1&=\alpha_1+\beta_1+h\gamma_1&=\;&\alpha_2+\beta_2+h\gamma_2,\\
k_2&=Q(\alpha_1,\beta_1,\gamma_1)&=\;&Q(\alpha_2,\beta_2,\gamma_2),\\
k_3&=C(\alpha_1,\beta_1,\gamma_1)&=\;&C(\alpha_2,\beta_2,\gamma_2).
\end{aligned}
\label{valkgen}
\end{equation}

In fact, we already know two solutions of Eqs.~\eqref{eqk1}, \eqref{eqk2} and \eqref{eqk3}, namely $(x,\,y,\,z)=(\alpha_1,\beta_1,\gamma_1)$ and  $(x,\,y,\,z)=(\alpha_2,\beta_2,\gamma_2)$.

To obtain a third solution of Eqs.~\eqref{eqk1}, \eqref{eqk2} and \eqref{eqk3}, we eliminate $x,\,y$ from these three equations when, in view of the relations~\eqref{condQC}, we get the following cubic equation in $x_3$:
\begin{equation}
(z-\gamma_1)(z-\gamma_2)(z-\gamma_3)=0, \label{cubiceqnz}
\end{equation}
where $\gamma_3$ is a rational function of the parameters occurring in the parametric solution of the simultaneous equations~\eqref{qdcubeqn} and \eqref{habgen}. The first two roots of Eq.~\eqref{cubiceqnz} namely, $z=\gamma_1$ and $z=\gamma_2$ yield the two known solutions of the simultaneous equations~\eqref{eqk1}, \eqref{eqk2} and \eqref{eqk3}.

 We will use the  third solution $z=\gamma_3$ to obtain  the diophantine  chains \eqref{qdchngen}, \eqref{cubchngen} and \eqref{chnhabc}. On substituting $z=\gamma_3$ in Eqs.~\eqref{eqk1} and \eqref{eqk2}, and eliminating $y$ from these two equations, we get a quadratic equation in $x$ which will have two rational roots if its discriminant is a perfect square. If we can choose the parameters such that  the discriminant is a perfect square, we will get two rational solutions of Eqs.~\eqref{eqk1} and \eqref{eqk2}, and we thus get a solution of the simultaneous diophantine chains~\eqref{qdchngen} and \eqref{cubchngen}, and hence also of the sextic   diophantine chain~\eqref{sexticchn}.

\subsection{}
When $\phi(x,\,y,\,z)$ is the sextic   form defined by \eqref{sexticphi}, we have the identity,
\begin{equation}
\phi(x,\,y,\,z)=(x^3+y^3-z^3)^2-4(xy)^3, \label{ident1}
\end{equation}
from which it follows that a solution of the simultaneous diophantine  chains,
\begin{align}
x_1^3+y_1^3-z_1^3 &=& x_2^3+y_2^3-z_2^3&=&x_3^3+y_3^3-z_3^3, \label{cubicchn}\\
x_1y_1&=&x_2y_2&=&x_3y_3, \label{qdchn}
\end{align}
will yield a solution of the sextic   diophantine chain~\eqref{sexticchn}

Following the method described in Section~2.1, we will first obtain a solution of the simultaneous diophantine  equations,
\begin{align}
x_1^3+y_1^3-z_1^3 &= x_2^3+y_2^3-z_2^3, \label{cubicab}\\
x_1y_1&=x_2y_2, \label{qdab}\\
x_1+y_1+hz_1&=x_2+y_2+hz_2, \label{hab}
\end{align}
where $h$ is a rational parameter.

The complete solution of  Eq.~\eqref{qdab} is given by 
\begin{equation}
x_1 = pu,\quad  y_1 = qv, \quad x_2 = pv, \quad y_2 = qu, \label{valab}
\end{equation}
where $p,\,q,\,u,\,v$ are arbitrary parameters.

With these values of $a_i,\,b_i,\;i=1,\,2$, Eq.~\eqref{cubicab} may be written as,
\begin{equation}
z_1^3-z_2^3-(u-v)(u^2+uv+v^2)(p-q)(p^2+pq+q^2)=0,\label{cubicab1}
\end{equation}
and on writing,
\begin{equation}
\begin{aligned}
z_1& = (n-m)(pu-qv)-m(pv+q(u+v)), \\
z_2 & = m(pu-qv)+n(pv+q(u+v)), \label{vala3b3} 
\end{aligned}
\end{equation}
where $m,\,n$ are  arbitrary parameters,
Eq.~\eqref{cubicab1} reduces to
\begin{multline}
(u^2+uv+v^2)(p^2+pq+q^2)\{(2m^3p+m^3q-3m^2np+3mn^2p\\
-n^3p+n^3q+p-q)u+(m^3p-m^3q+3m^2nq-3mn^2q\\
+n^3p+2n^3q-p+q)v\}=0.
\end{multline}
Accordingly, we get,
\begin{equation}
\begin{aligned}
u&=m^3p-m^3q+3m^2nq-3mn^2q+n^3p+2n^3q-p+q,\\
v&=-(2m^3p+m^3q-3m^2np+3mn^2p-n^3p+n^3q+p-q),
\end{aligned} \label{valuv}
\end{equation}
and on substituting these values of $u$ and $v$ in \eqref{valab} and \eqref{vala3b3}, we get a solution of the simultaneous equations~\eqref{cubicab} and \eqref{qdab} which may be written in terms of arbitrary parameters $m,\,n,\,p$ and $q$ as  $(x_i,\,y_i,\,z_i)=(\alpha_i,\,\beta_i,\,\gamma_i),\;i=1,\,2$, where 
\begin{equation}
\begin{aligned}
\alpha_1 &= (m^3+n^3-1)p^2+(-m^3+3m^2n-3mn^2+2n^3+1)pq,\\
 \beta_1 &= (-2m^3+3m^2n-3mn^2+n^3-1)pq+(-m^3-n^3+1)q^2,\\
 \gamma_1 &= (m^4-2m^3n+3m^2n^2-2mn^3+n^4+2m-n)p^2\\
& \quad \quad +(m^4-2m^3n+3m^2n^2-2mn^3+n^4-m+2n)pq\\
& \quad \quad+(m^4-2m^3n+3m^2n^2-2mn^3+n^4-m-n)q^2, \\
\alpha_2 &= (-2m^3+3m^2n-3mn^2+n^3-1)p^2+(-m^3-n^3+1)pq, \\
\beta_2 &= (m^3+n^3-1)pq+(-m^3+3m^2n-3mn^2+2n^3+1)q^2, \\
\gamma_2 &= (m^4-2m^3n+3m^2n^2-2mn^3+n^4-m-n)p^2\\
& \quad \quad+(m^4-2m^3n+3m^2n^2-2mn^3+n^4+2m-n)pq\\
& \quad \quad+(m^4-2m^3n+3m^2n^2-2mn^3+n^4-m+2n)q^2.
\end{aligned} \label{solab}
\end{equation}

We note that the solution~\eqref{solab} also satisfies Eq.~\eqref{hab} when $h=-(m^2-mn+n^2)$.

We will now obtain a solution of the simultaneous diophantine  chains \eqref{cubicchn}, \eqref{qdchn}  and the auxiliary diophantine    chain~\eqref{chnhabc} by obtaining 
three distinct solutions of the simultaneous diophantine  equations~\eqref{eqk1}, \eqref{eqk2}, \eqref{eqk3}
where $h=-(m^2-mn+n^2),\;k_1=\alpha_1+\beta_1+h\gamma_1,\;k_2=\alpha_1\beta_1,\;k_3=\alpha_1^3+\beta_1^3-\gamma_1^3$, with $\alpha_1,\;\beta_1,\;\gamma_1$ being defined  by \eqref{solab}, so that Eqs.~\eqref{eqk1}, \eqref{eqk2}, \eqref{eqk3}  have already two known solutions $(x,\,y,\,z)=(\alpha_1,\,\beta_1,\,\gamma_1)$ and $(x,\,y,\,z)=(\alpha_2,\,\beta_2,\,\gamma_2)$.

To obtain a third solution of the simultaneous diophantine  equations \eqref{eqk1}, \eqref{eqk2}, \eqref{eqk3}, we eliminate $x$ and $y$ from these three equations to get the  cubic equation \eqref{cubiceqnz} where $\gamma_1,\;\gamma_2$ are defined by \eqref{solab} and 
\begin{multline}
\gamma_3=\{(m^{10}-5m^9n+15m^8n^2-30m^7n^3+45m^6n^4-51m^5n^5+45m^4n^6\\
-30m^3n^7+15m^2n^8-5mn^9+n^{10}+2m^7-10m^6n+24m^5n^2-38m^4n^3\\
+40m^3n^4-30m^2n^5+14mn^6-4n^7+5m^4-{10}m^3n+15m^2n^2-10mn^3\\
+5n^4+m-2n)p^2+(m^{10}-5m^9n+15m^8n^2-30m^7n^3+45m^6n^4-51m^5n^5\\
+45m^4n^6-30m^3n^7+15m^2n^8-5mn^9+n^{10}+5m^7-19m^6n+42m^5n^2\\
-59m^4n^3+58m^3n^4-39m^2n^5+17mn^6-4n^7+2m^4-4m^3n+6m^2n^2\\
-4mn^3+2n^4+m+n)pq+(m^2-mn+n^2-1)(m^4-2m^3n+3m^2n^2-2mn^3\\
+n^4+m^2-mn+n^2+1)(m^4-2m^3n+3m^2n^2-2mn^3+n^4+2m-n)q^2\}\\
\times \{(m^2-mn+n^2-1)((m^2-mn+n^2)^2+m^2-mn+n^2+1)\}^{-1}. \label{valc3}
\end{multline}

While the first two roots, $z=\gamma_1$ and $z=\gamma_2$, of Eq.~\eqref{cubiceqnz} yield the two known solutions of Eqs.~\eqref{eqk1}, \eqref{eqk2}, \eqref{eqk3}, the third root $z=\gamma_3$ will yield a new solution. We will take  $z=\gamma_3$ in Eqs.~\eqref{eqk1} and \eqref{eqk2}, and solve them to get the values of $x$ and $y$. It  follows from  \eqref{eqk1} that $x+y=k_1-h\gamma_3$, and hence $(x-y)^2=(x+y)^2-4xy=(k_1-h\gamma_3)^2-4k_2$. Thus, both $x,\,y$ will be rational if $(k_1-h\gamma_3)^2-4k_2$ is a perfect square. This gives us a quartic function in $p$ and $q$ to be made a perfect square. As this function is too cumbersome to write in full, we restrict ourselves to writing it as follows:
\begin{multline}
\{(m-2n)(m^2-mn+n^2)^4+5(m^2-mn+n^2)^3+2(m-2n)(m^2-mn+n^2)+1\}^2p^4\\
+\cdots+(m^2-mn+n^2-1)^2(2m^3-3m^2n+3mn^2-n^3+1)^2\\
\times \{(m^2-mn+n^2)^2+m^2-mn+n^2+1\}^2q^4. \label{quarticpq}
\end{multline}

Since the coefficients of $p^4$ and $q^4$ in the quartic function~\eqref{quarticpq} are perfect squares, we can readily find infinitely many values of $p$ and $q$ that make the function~\eqref{quarticpq} a perfect square by repeatedly applying  a method described by Fermat (as quoted by Dickson \cite[p. 639]{Di}), and each such solution will lead to a solution of the simultaneous diophantine chains \eqref{cubicchn}, \eqref{qdchn}  and \eqref{chnhabc}, and hence also of the sextic  chain~\eqref{sexticchn}, in terms of the rational parameters $m$ and $n$.  

 As an example, the quartic function~\eqref{quarticpq} becomes a perfect square  if we choose $p$ and $q$ as follows:
\begin{equation}
\begin{aligned}
p& = -(m-n)((m+n)t+2)(t^3-1),\\
 q& =((2m-n)t+1)((m-n)t^3+t^2+m),
\end{aligned}
\label{valpq}
\end{equation}
where 
\begin{equation}
t=m^2-mn+n^2 \label{valt}. 
\end{equation}

This yields a solution of the sextic  chain~\eqref{sexticchn} which may be written, in terms of arbitrary parameters $m,\,n$, as follows:
\begin{equation}
\begin{aligned}
x_1&=(m-n)((m+n)t+2)(t^3-1)(3(m-n)t^6+(2m-n)(m-2n)t^4\\
& \quad \quad -3(m^3+mn^2-n^3)t^2-3tm^2+m-2n),\\
y_1&=((m-2n)t-1)((2m-n)t+1)^2(2t^2+m-2n)((m-n)t^3+t^2+m),\\
z_1&=3(m-n)^2t^{11}+9m(m-n)^2t^9+(19m^4-68m^3n\\ & \quad \quad
+81m^2n^2-35mn^3+4n^4)t^7+(4m^5-67m^4n+163m^3n^2\\ & \quad \quad
-176m^2n^3+104mn^4-26n^5)t^5-(23m^4-22m^3n-18m^2n^2\\ & \quad \quad
+32mn^3-13n^4)t^4-(19m^3-36m^2n+39mn^2-14n^3)t^3\\ & \quad \quad
+t(2m^4-7m^3n-6m^2n^2+8mn^3-4n^4)\\
& \quad \quad+(m-2n)(5m^2-5mn+2n^2),\\
x_2&=-(m-n)((m+n)t+2)((m-2n)t-1)((2m-n)t+1)\\
& \quad \quad \times (2t^2+m-2n)(t^3-1),\\
y_2&=-((2m-n)t+1)((m-n)t^3+t^2+m)(3(m-n)t^6\\
& \quad \quad+(2m-n)(m-2n)t^4-3(m^3+mn^2-n^3)t^2-3m^2t+m-2n),\\
z_2&=3(m-n)^2t^{11}-(14m^4-37m^3n+36m^2n^2-22mn^3+8n^4)t^7\\ & \quad \quad
-(23m^5-98m^4n+140m^3n^2-103m^2n^3+37mn^4-4n^5)t^5\\ & \quad \quad
+(m^4+58m^3n-90m^2n^2+46mn^3-5n^4)t^4+(23m^3-12m^2n\\ & \quad \quad
-9mn^2+8n^3)t^3+(11m^4-28m^3n+39m^2n^2-19mn^3+2n^4)t\\ & \quad \quad
-(m-2n)(m^2+2mn-2n^2),\\
x_3&=((m-2n)t-1)((m-n)t^3+t^2+m)(3(m-n)t^6\\
& \quad \quad+(2m-n)(m-2n)t^4-(3(m^3+mn^2-n^3))t^2-3m^2t+m-2n),\\
y_3&=(m-n)((m+n)t+2)((2m-n)t+1)^2(2t^2+m-2n)(t^3-1),\\
z_3&=3(m-n)^2t^{11}+9(m-n)^3t^9+(13m^4-35m^3n+36m^2n^2\\ & \quad \quad
-14mn^3+n^4)t^7+(19m^5-55m^4n+79m^3n^2-44m^2n^3\\ & \quad \quad
-mn^4+7n^5)t^5+(4m^4-17m^3n+54m^2n^2-41mn^3+10n^4)t^4\\& \quad \quad
-2(11m^3-21m^2n+5n^3)t^3 -t(22m^4-68m^3n+78m^2n^2\\& \quad \quad
-47mn^3+10n^4)-(m-2n)(4m^2-7mn+4n^2),
\end{aligned}
\label{sexticchnsol1}
\end{equation}
where, as before, $t=m^2-mn+n^2$.

As a numerical example, when $m=1,\,n=2$, after removing common factors, we get the following solution of the sextic  chain~\eqref{sexticchn}:
\begin{equation}
\begin{aligned}
x_1&=100958,\quad &y_1&=425,\quad &z_1&=113259,\\
x_2&= -7150,\quad &y_2&= -6001,\quad &z_2&= 75081, \\
x_3&=-60010,\quad &y_3&= -715,\quad &z_3&=59223. 
\end{aligned}
\label{sexticchnsol1ex1}
\end{equation}

\subsection{}
It is interesting to note that when $\phi(x,\,y,\,z)$ is defined by \eqref{sexticphi}, in addition to the identity \eqref{ident1}, we also have  the identity,
\begin{equation}
\phi(x,\,y,\,z)=4Q^3(x,\,y,\,z)-3C^2(x,\,y,\,z),
\end{equation}
where
\begin{equation}
\begin{aligned}
Q(x,\,y,\,z)&=x^2+y^2+z^2+xy+yz+zx,\\
C(x,\,y,\,z)&=x^3+y^3+z^3+2x^2y+2xy^2+2x^2z+2xz^2\\
& \quad \quad +2y^2z+2yz^2+2xyz.
\end{aligned}
\end{equation}

We can now apply the method described in Section~2.1 to obtain solutions of the diophantine chain~\eqref{sexticchn}.

A parametric solution of the simultaneous diophantine  equations~\eqref{qdcubeqn}  may be obtained by a straightforward application of  the method described in \cite{Ch0}. We accordingly omit the tedious details and simply state below the  solution thus obtained.

If we define three functions $f_i(u,\,v,\,w),\; i=1,\,2,\,3$, as 
\begin{equation}
\begin{aligned}
f_1(u,\,v,\,w)&=(3u^2-2uv-2uw-v^2+2vw-w^2)(u^3+uv^2\\
& \quad \quad -2uvw+uw^2-2v^3+2v^2w+2vw^2-2w^3),\\
f_2(u,\,v,\,w)&=-2(v-w)(u+v-w)(u-v-w)(u-v+w)\\
& \quad \quad \times(uv+uw-v^2-vw-w^2),\\
f_3(u,\,v,\,w)&=u^6-2u^5v-2u^5w+2u^4v^2+2u^4w^2-2u^3v^3\\
& \quad \quad +2u^3v^2w+2u^3vw^2-2u^3w^3+2u^2v^4+2u^2v^3w\\
& \quad \quad -8u^2v^2w^2+2u^2vw^3+2u^2w^4-2uv^5+2uv^3w^2\\
& \quad \quad +2uv^2w^3-2uw^5+v^6-2v^5w+2v^4w^2-2v^3w^3\\
& \quad \quad +2v^2w^4-2vw^5+w^6,
\end{aligned}
\end{equation}
the aforesaid solution of the simultaneous equations~\eqref{qdcubeqn} may be written, in terms of arbitrary parameters $p,\,q,\,r$ and $m$ as $(x_i,\,y_i,\,z_i)=(\alpha_i,\,\beta_i,\,\gamma_i),\;i=1,\,2$, where
\begin{equation}
\begin{aligned}
\alpha_1&=f_1(p,\,q,\,r)m^2+f_2(p,\,q,\,r)m+pf_3(p,\,q,\,r),\\
\beta_1&=f_1(q,\,r,\,p)m^2+f_2(q,\,r,\,p)m+qf_3(p,\,q,\,r),\\
\gamma_1&=f_1(r,\,p,\,q)m^2+f_2(r,\,p,\,q)m+rf_3(p,\,q,\,r),\\
\alpha_2&=f_1(p,\,q,\,r)m^2-f_2(p,\,q,\,r)m+pf_3(p,\,q,\,r),\\
\beta_2&=f_1(q,\,r,\,p)m^2-f_2(q,\,r,\,p)m+qf_3(p,\,q,\,r),\\
\gamma_2&=f_1(r,\,p,\,q)m^2-f_2(r,\,p,\,q)m+rf_3(p,\,q,\,r).
\end{aligned}
\label{sol2xyz}
\end{equation}

Further, the values of $\alpha_i,\,\beta_i,\,\gamma_i,\;i=1,\,2$, given by \eqref{sol2xyz} satisfy Eq.~\eqref{habgen} where
\begin{equation}
h=(p^2+q^2-r^2)/(p^2+pq-pr+q^2-qr). \label{sol2valh}
\end{equation}

Now with the values of $\alpha_i,\,\beta_i,\,\gamma_i,\;i=1,\,2$ and $h$ given by \eqref{sol2xyz} and \eqref{sol2valh}, we will solve Eqs.~\eqref{eqk1}, \eqref{eqk2} and \eqref{eqk3}. On eliminating $x$ and $y$ from these three equations, we get Eq.~\eqref{cubiceqnz} where 
\begin{multline}
\gamma_3=(p^2+pq-pr+q^2-qr)\{(3p^7-2p^6q-2p^6r+4p^5q^2-4p^5qr\\
-5p^4q^3+2p^4q^2r+6p^4qr^2-3p^4r^3-5p^3q^4+8p^3q^3r-6p^3q^2r^2\\
+8p^3qr^3-5p^3r^4+4p^2q^5+2p^2q^4r-6p^2q^3r^2-10p^2q^2r^3\\
+10p^2qr^4-2pq^6-4pq^5r+6pq^4r^2+8pq^3r^3+10pq^2r^4-36pqr^5\\
+18pr^6+3q^7-2q^6r-3q^4r^3-5q^3r^4+18qr^6-11r^7)m^2\\
+(p^3+q^3-r^3)f_3(p,\,q,\,r)\}\{p^4+p^3q-p^3r+3p^2q^2\\
-3p^2qr+pq^3-3pq^2r+3pqr^2-pr^3+q^4-q^3r-qr^3+r^4\}^{-1}.
\end{multline}

With $z=\gamma_3$, we have to solve Eqs.~\eqref{eqk1} and \eqref{eqk2}. This leads to a quadratic equation in $x,\,y$ whose discriminant is to be made a perfect square. As this discriminant is too cumbersome to write, we will take specific numerical values of the parameters that yield the desired sextic   diophantine chains~\eqref{sexticchn}. We take for simplicity $q=0$, and now the condition that the discriminant be a perfect square reduces to finding rational solutions of the following equation:
\begin{multline}
Y^2=r(36p^{11}+96p^{10}r+220p^9r^2+357p^8r^3+522p^7r^4+541p^6r^5\\
\;\;+462p^5r^6+228p^4r^7+22p^3r^8-99p^2r^9-54pr^{10}-27r^{11})m^4\\
\;\;\;+2(8p^{10}+28p^9r+12p^8r^2+12p^7r^3-9p^6r^4+12p^5r^5-8p^4r^6\\
+14p^3r^7-12p^2r^8-9r^{10})(p^2+pr+r^2)^2m^2\\
+r(4p^3-3r^3)(p^2-pr+r^2)^2(p^2+pr+r^2)^4. \label{ecsol2}
\end{multline}

When  $p=3,\,r=4$, Eq.~\eqref{ecsol2} reduces to the quartic equation,
\begin{equation}
Y^2=4916053296m^4-16574603472m^2-106422358224, \label{ecsol2ex1}
\end{equation}
which represents a quartic model of an elliptic curve. On making the birational transformation defined by the relations,
\begin{equation}
 \begin{aligned}
m &= (37/3)(521u-11v+318899904)/(42871u-11v-2751064896),\\
Y& = 958300(1331u^3-289453824u^2-68536946349036u+86313500544v\\
& \quad \quad +3183632918644552704)/(42871u-11v-2751064896)^2
\end{aligned}
\label{biratmY}
\end{equation}
and 
\begin{equation}
\begin{aligned}
u &= (136052568m^2+1925Y-24886644m-96169512)/\{2(3m-37)^2\},\\
v& = 175(9012764376m^3+128613Ym-1231896204m^2-19277Y\\
& \quad \quad -15193386516m-15819539736)/\{2(3m-37)^3\},
\end{aligned}
\label{biratuv}
\end{equation}
Eq.~\eqref{ecsol2ex1} reduces to the Weierstrass form of the elliptic curve given by
\begin{equation}
v^2 = u^3+u^2+51492677220u-3062315437673472. \label{ecsol2ex1w}
\end{equation}
The rank of the elliptic curve~\eqref{ecsol2ex1w}, as determined by the software {\sf {SAGE}}, is  3, with the three generators of the Mordell-Weil group being
\[
\begin{aligned}
&(101376,\, 56565600),\quad (3761676,\, -7308840000),\\
&\mbox{\rm and} \quad (498157004/529,\, 11417003301600/12167).
\end{aligned}
\] 
 We can now find infinitely many rational points on the elliptic curve~\eqref{ecsol2ex1w} using the group law, and then find infinitely many rational points on the quartic curve~\eqref{ecsol2ex1} using the relations~\eqref{biratmY}. These rational points on the curve~\eqref{ecsol2ex1w} yield infinitely many numerical examples of the sextic  chain~\eqref{sexticchn}.

In order to find rational points of small height on the curve~\eqref{ecsol2ex1}, we used Stoll's program `ratpoints' \cite{St} and readily  obtained the following four  values of $m$ for which the right-hand side of Eq.~\eqref{ecsol2ex1} becomes a perfect square:
\[37/3, \quad 481/87, \quad 14911/4695, \quad 135679/50151.\]
The first two values of $m$ do not lead to nontrivial sextic   chains but the next two values yield the following two solutions of the sextic  chain~\eqref{sexticchn}:
\begin{equation}
\begin{aligned}
(x_1,\,y_1,\,z_1)&=(14900543,\,-2461462,\,15194895),\\
(x_2,\,y_2,\,z_2)&=(12571823,\,2923703,\,13884990),\\
(x_3,\,y_3,\,z_3)&=(4528874,\,11547071, \,13636239),
\end{aligned}
\label{sexticchnsol2ex1}
\end{equation}
and
\begin{equation}
\begin{aligned}
(x_1,\,y_1,\,z_1)&=(17217348683,\, -3153451318,\, 17759190363),\\
  (x_2,\,y_2,\,z_2)&=(14274889211,\, 3650986211,\,16104056910),\\
(x_3,\,y_3,\,z_3)&=(11570059211,\,  7442013386,\,  15638543835).
\end{aligned}
\label{sexticchnsol2ex2}
\end{equation}

\section{Mordell curves related to sextic   diophantine  chains}
With every solution of the sextic   diophantine chain~\eqref{sexticchn}, we may associate a Mordell curve~\eqref{Mordelleq} where we take 
\begin{equation}
4k=\phi(x_1,\,y_1,\,z_1)=\phi(x_2,\,y_2,\,z_2)=\phi(x_3,\,y_3,\,z_3), \label{relk}
\end{equation}
and, as noted in the Introduction, there will, in general, be 9 rational  points on the curve~\eqref{Mordelleq} whose coordinates are immediately obtained from the sextic   diophantine  chain. 

We now consider the family of Mordell curves related to the parametric solution~\eqref{sexticchnsol1} of the sextic   diophantine chain~\eqref{sexticchn}. Using the relations~\eqref{sexticchnsol1} and \eqref{relk}, we can compute the value of $k$ in terms of the arbitrary parameters $m$ and $n$. We thus get a  Mordell curve defined  over the field $\mathbb{Q}(m,\,n)$ on which we can readily find 9 rational points. The value of $k$ is too cumbersome to write and we do not give it explicitly.  

We note that in view of the relations~\eqref{qdab}, only 7 of the 9 known rational points on our Mordell curve    are actually distinct. We will now  apply a theorem of Silverman \cite[Theorem 11.4, p. 271]{Sil} to show that 6 of these 7 points are linearly independent in the group of rational points of the Mordell curve. For this,  we must find specific numerical values of $m$ and $n$ such that 6 of the 7 points are linearly independent on the related Mordell curve  over $\mathbb{Q}$.

When $m=1,\,n=2$, the numerical solution~\eqref{sexticchnsol1ex1} of the sextic  chain~\eqref{sexticchn} is related to the Mordell curve,
\begin{equation}
y^2=x^3+44906825622115054978352852841, \label{Mordelleqex1}
\end{equation}
on which we get 7 rational points whose coordinates are given below:
\begin{align*}
&(42907150, -211912492824721), \quad &(48135075, 211912569590346),\\
 &(11434402122, 1240928701242633), \quad &(-450561081, 211696384806720),\\
 &(-536829150, 211546966949721), \quad &(-42344445, 211912127298846), \\
&(-3553972230, -4195525176279). &
\end{align*}
The regulator of the first six of these points, as computed by {\sf {SAGE}}, is 10390179.16.  As this is nonzero,  it follows from a well-known theorem \cite[Theorem 8.1, p. 242]{SZ} that these 6 points are linearly independent.  It follows that the generic rank of the family of Mordell curves related to the sextic   chain given by the parametric solution~\eqref{sexticchnsol1} is at least 6.

We could not find any numerical values of $m$ and $n$ such that all the 7 known points are linearly independent in the group of rational points of the Mordell curve. 

Next we consider the Mordell curves related to the two numerical solutions~\eqref{sexticchnsol2ex1} and \eqref{sexticchnsol2ex2} of the sextic  chain~\eqref{sexticchn}.

The Mordell curve related to the solution  \eqref{sexticchnsol2ex1} is
\begin{equation}
y^2=x^3+60881141602872940726223731917150516833400, \label{Mordelleqex2}
\end{equation}
on which we get  9 distinct rational points $P_1,\,P_2,\,\ldots,\,P_9$ whose coordinates are as follows:
\begin{align*}
&P_1=(-36677120373866, 107436818637424863748), \\& P_2=(226412186327985, -3415747486107335266755), \\&P_3= (-37401656636490, -92523324542363200620), \\& P_4=(36756276620569, 332475185129096665153), \\&P_5= (174559636636770, -2319461032255991683920), \\&P_6= (40595586917970, -357467113076423815080), \\&P_7= (52295229628054, 451550293014200834692), \\&P_8= (61756808264886, -544440667643008046316), \\&P_9= (157458619905969, -1991177257485343473603).
 \end{align*}
The regulator of the 6 points $P_1,\,P_2,\,P_3,\,P_4,\,P_5$ and $P_8$ is 11390832.16. Since this is nonzero, these 6 points are independent, and the rank of the Mordell curve~\eqref{Mordelleqex2} is at least 6.

Similarly, the solution~\eqref{sexticchnsol2ex2} of the sextic  chain~\eqref{sexticchn}  yields 9 distinct rational points on  the 
Mordell curve, 
\begin{equation}
\begin{aligned}
y^2&=x^3+292994052252973481957957313165581343604439577105\\
& \quad \quad 2383353330586891185336.
\end{aligned}
\label{Mordelleqex3}
\end{equation}
Again we found that only 6 of the 9 points are independent, and so the rank of the Mordell curve~\eqref{Mordelleqex3}
is at least 6.

We could not determine the precise rank of any of the three Mordell curves \eqref{Mordelleqex1}, \eqref{Mordelleqex2} or \eqref{Mordelleqex3} as the value of $k$ for each of these three curves is very large. 

\section{An open problem}
It would be of interest to solve the sextic   diophantine  chain,
\begin{equation}
\phi(x_1,\,y_1,\,z_1)=\phi(x_2,\,y_2,\,z_2)=\cdots=\phi(x_n,\,y_n,\,z_n), \label{sexticchngen}
\end{equation}
when $\phi(x,\,y,\,z)$ is defined by \eqref{sexticphi} and $n > 3$. While the existence of such diophantine  chains when $n=4$ and $n=5$ is not inconceivable, it certainly seems that there must be an upper bound for $n$ for the solvability of the diophantine  chain  \eqref{sexticchngen}. It would be of interest to determine the largest integer $n$ for which the diophantine chain~\eqref{sexticchngen} is solvable.

Any solution of the diophantine chain~\eqref{sexticchngen} will immediately yield $3n$ rational points on the Mordell curve~\eqref{Mordelleq} where $k=\phi(x_1,\,y_1,\,z_1)/4$, and may therefore yield Mordell curves of rank  higher than the examples already known in the literature.

\begin{center}
\Large
Acknowledgment
\end{center}

The first author thanks the Harish-Chandra Research Institute, Prayagraj for providing him with all necessary facilities that have helped him to pursue his research work in mathematics.

\medskip

\noindent Postal Address: Ajai Choudhry, 
\newline \hspace{1.05 in}
13/4 A Clay Square,
\newline \hspace{1.05 in} Lucknow - 226001, INDIA.
\newline \noindent  E-mail: ajaic203@yahoo.com
\vspace{.2cm} \\
\noindent Postal Address: Arman Shamsi Zargar, 
\newline \hspace{1.05 in}
Department of Mathematics and Applications,
\newline \hspace{1.05 in} 
Faculty of Science,
\newline \hspace{1.05 in} 
University of Mohaghegh Ardabili,
\newline \hspace{1.05 in} 
Ardabil 56199-11367, IRAN.
\newline \noindent  E-mail: zargar@uma.ac.ir

\end{document}